\documentclass[11pt,leqno]{amsart}

\usepackage{amscd}
\usepackage{color}
\usepackage{amssymb}
\usepackage{amsfonts}
\usepackage{latexsym}
\usepackage{verbatim}

\theoremstyle{plain}
\newtheorem{theorem}{Theorem}[section]
\newtheorem{definition}[theorem]{Definition}

\renewcommand{\b}{\begin{equation}}
\newcommand{\e}{\end{equation}}

\thanks{This work was supported by GNSAGA of INdAM}
\address{Dipartimento di Ingegneria e Scienze dell'Informazione e Matematica \\ Universit\`a dell'Aquila\\
via Vetoio\\ 67100 L'Aquila\\ Italy}
\email{lucio.bedulli@univaq.it}
\address{Dipartimento di Matematica G. Peano \\ Universit\`a di Torino\\
Via Carlo Alberto 10\\
10123 Torino\\ Italy}
 \email{luigi.vezzoni@unito.it}

\author{Lucio Bedulli and Luigi Vezzoni}
\date{\today}

\title{A remark on the Laplacian flow and the modified Laplacian co-flow in ${\rm G}_2$-Geometry}
\begin{document}

\maketitle
\begin{abstract}
We observe that the DeTurck Laplacian flow of ${\rm G}_2$-structures introduced by Bryant and Xu as a gauge fixing of the Laplacian flow can be regarded as 
a flow of ${\rm G}_2$-structures (not necessarily closed) which fits in the general framework introduced by Hamilton in \cite{positive}.  
\end{abstract}

\section{Introduction}
In \cite{Bryant} Bryant introduced a geometric flow in ${\rm G}_2$-geometry which evolves an initial closed ${\rm G}_2$-structure $\varphi_0$ in the direction of its Laplacian. 

Given  a compact $7$-dimensional manifold with a closed ${\rm G}_2$-structure $(M,\varphi_0)$, a {\em Laplacian flow} is a solution to the evolution equation 
\begin{equation}\label{LF}
\tfrac{\partial}{\partial t}\varphi_t=\Delta_{\varphi_t}\varphi_t\,,\quad d\varphi_t=0\,,\quad  \varphi_{|t=0}=\varphi_0\,.
\end{equation}
The well-posedness of equation \eqref{LF} is proved in \cite{BryantXu} by applying the Nash-Moser theorem to the gauge fixing 
\begin{equation}\label{DTLF}
\tfrac{\partial}{\partial t}\varphi_t=\Delta_{\varphi_t}\varphi_t+\mathcal{L}_{V(\varphi_t)}\varphi_t\,,\quad d\varphi_t=0\,,\quad  \varphi_{|t=0}=\varphi_0\,,
\end{equation}
where $\mathcal L$ is the Lie derivative and $V\colon C^{\infty}(M,\Lambda^3_+)\to C^{\infty}(M,TM)$ is a first order
differential operator which depends on the choice of a connection on $M$. Here $\Lambda^3_+$ denotes the open subbundle of $\Lambda^3$ of ${\rm G}_2$-structures on $M$.
A solution to \eqref{DTLF} is usually called a {\em DeTurck Laplacian flow}. 

A  DeTurck Laplacian flow $\varphi_t$ is also a solution to  
\begin{equation}\label{flow}
\tfrac{\partial}{\partial t}\varphi_t=dd^*_{\varphi_t}\varphi_t+d\iota_{V(\varphi_t)}\varphi_t\,,\quad  \varphi_{|t=0}=\varphi_0\,.
\end{equation}

In the present note we observe that equation \eqref{flow} fits in the general framework introduced by Hamilton in \cite{positive}. As a direct consequence we have the following  theorem which in particular implies the well-posedness of 
\eqref{DTLF} 
\begin{theorem}\label{main1}
Let $(M,\varphi_0)$ be a compact $7$-dimensional manifold with a ${\rm G}_2$-structure. Then equation \eqref{flow} has a unique short-time solution.
\end{theorem}
 
In \cite{K} Karigiannis, McKay and Tsui introduced the {\em Laplacian co-flow} as the solution to the evolution equation 
\begin{equation}\label{LCF}
\tfrac{\partial}{\partial t}(*_{\varphi_t}\varphi_t)=-\Delta_{\varphi_t}*_{\varphi_t}\varphi_t\,,\quad d*_{\varphi_t}\varphi_t
=0\,,\quad  \varphi_{|t=0}=\varphi_0\,.
\end{equation}
where in this case $\varphi_0$ is supposed to be co-closed with respect to the metric induced by itself. 
The well-posedness of this last equation is still an open problem and Grigorian introduced in \cite{Gri} the following modification
\begin{equation}\label{MLCF}
\tfrac{\partial}{\partial t}(*_{\varphi_t}\varphi_t)= \Delta_{\varphi_t}*_{\varphi_t}\varphi_t+2d((A-{\rm Tr}(T(\varphi_t))\varphi_t )\,,\quad d*_{\varphi_t}\varphi_t
=0\,,\quad  \varphi_{|t=0}=\varphi_0\,,
\end{equation}
where $A$ is a constant and $T(\varphi_t)$ is the torsion of $\varphi_t$. In \cite{Gri} it is proved the well-posedness of \eqref{MLCF} following the same approach of Bryant in \cite{Bryant} by applying the Nash-Moser theorem  to the gauge fixing 
\begin{equation}\label{DTMLCF}
\tfrac{\partial}{\partial t}(*_{\varphi_t}\varphi_t)=\Delta_{\varphi_t}*_{\varphi_t}\varphi_t+2d((A-{\rm Tr}(T(\varphi_t))\varphi_t )+\mathcal{L}_{V(\varphi_t)}\varphi_t\,,\quad d*_{\varphi_t}\varphi_t=0\,,\quad  \varphi_{|t=0}=\varphi_0\,,
\end{equation} 
Any solution to this last equation \eqref{DTMLCF}
satsfies
\begin{equation}\label{DTMLCF}
\tfrac{\partial}{\partial t}(*_{\varphi_t}\varphi_t)=dd^{*}_{\varphi_t}*_{\varphi_t}\varphi_t+2d((A-{\rm Tr}(T(\varphi_t))\varphi_t )+d\iota_{V(\varphi_t)}\varphi_t\,,\quad  \varphi_{|t=0}=\varphi_0\,,
\end{equation}
Analogously to theorem \ref{main1} we have 
\begin{theorem}\label{main2}
Let $(M,\varphi_0)$ be a compact $7$-dimensional manifold with a ${\rm G}_2$-structure. Then equation \eqref{DTMLCF} has a unique short-time solution.
\end{theorem}

\section{Proof of the results}
Both theorems \ref{main1} and \ref{main2} can be proved by using the following set-up introduced by Hamilton in \cite{positive}.

\medskip 
Let $M$ be an oriented compact manifold, $F$ a vector bundle over $M$, $U$ an open subbundle of $F$ and  
$$
E\colon C^{\infty}(M,U)\to  C^{\infty}(M,F)
$$
a second order differential operator. For $f\in C^{\infty}(M,U)$, we denote by $D E(f)\colon C^{\infty}(M,F)\to C^{\infty}(M,F)$ the linearization of $E$ at $f$ and by 
$\sigma D E(f)$ the principal symbol of
$D E(f)$.

\begin{definition}
{\em An integrability condition for $E$ is a first order linear differential operator 
$$ 
L\colon C^{\infty}(M,F)\to C^{\infty}(M,G)\,,
$$
where $G$ is another vector bundle over $M$, such that $L(E(f))=0$ for all $f\in C^{\infty}(M,U)$, and all the eigenvalues of $\sigma DE(f)$ restricted to $\ker \sigma L$
 have strictly positive real part. }
\end{definition}

\begin{theorem}[Hamilton {\cite[Theorem 5.1]{positive}}]\label{Ham_int}
Assume that $E$ admits
 an integrability condition.  Then for every $f_0\in C^{\infty}(M,U)$ the geometric flow 
\begin{equation}\label{flow_Ham}
\frac{\partial f}{\partial t}=E(f)\,,\quad f(0)=f_0\,,
\end{equation}
has a unique short-time solution. 
\end{theorem}

Now we can focus on the set-up of theorem \ref{main1}. Here we consider
$$
F=\Lambda^3\,,\quad U=\Lambda^3_+\,,\quad G=\Lambda^4\,,\quad E(\varphi)= dd^*_{\varphi}\varphi+d\iota_{V(\varphi)}\varphi\,,\quad 
L=d\colon C^{\infty}(M,\Lambda^3)\to C^{\infty}(M,\Lambda^4)\,.
$$
From \cite{BryantXu} it follows that for every $\varphi\in C^{\infty}(M,U)$ and every {\em closed} $\psi \in C^{\infty}(M,\Lambda^{3})$  we have 
$$
DE(\varphi)(\psi)=-\Delta_{\varphi} \psi+{\rm l.o.t.} 
$$
Hence all the assumptions of Hamilton's theorem \ref{Ham_int} are satisfied and theorem \ref{main1} follows. 

\medskip 
Notice that if the starting form $\varphi_0$ is closed, then the solution to \eqref{flow} is closed for every $t$ since 
$$
d\tfrac{\partial}{\partial t}\varphi=0\,.
$$
Therefore if $\varphi_0$ is closed, 
the unique solution $\varphi_t$ to \eqref{flow} solves also the DeTurck-Laplacian flow \eqref{DTLF}
and the short-time existence of the DeTurck-Laplacian flow \eqref{DTLF} can be deduced from Theorem \ref{main1}.

\bigskip 
About the proof of theorem \ref{main2} we set   

$$
\begin{aligned}
& F=\Lambda^4\,,\quad U=\Lambda^4_+\,,\quad G=\Lambda^5\,,\quad E(*_\varphi\varphi)= dd^*_{\varphi}\varphi+d\iota_{V(\varphi)}+2d((A-{\rm Tr}(T(\varphi))\varphi)\,,\\ 
&L=d\colon C^{\infty}(M,\Lambda^4)\to C^{\infty}(M,\Lambda^5)\,.
\end{aligned}
$$
From \cite{Gri} it follows 
$$
DE(*_\varphi\varphi)(\psi)=-\Delta_{\varphi} \psi+{\rm l.o.t.} 
$$
for every closed $\psi\in C^\infty(M,\Lambda^4)$ and the proof of theorem \ref{main2} follows. 

\bigskip
\noindent {\bf Acknowledgments.} We would like to thank Jason Lotay for useful conversations.

\end{document}